\documentclass[numbers,webpdf,imanum]{mod-ima-authoring-template}%
\usepackage{cases}
\usepackage{algorithm}
\usepackage{tikz}
\usetikzlibrary{calc,patterns,decorations.pathmorphing,decorations.markings}

\graphicspath{{Fig/}}
\theoremstyle{thmstyletwo}%
\newtheorem{theorem}{Theorem}
\newtheorem{lemma}[theorem]{Lemma}%

\numberwithin{equation}{section}

\def\nbyn{n\times n}
\def\Dext{D_\mathrm{ext}}
\def\Dint{D_\mathrm{int}}

\def\res{\mathrm{res}}
\def\tol{\mathrm{tol}}
\def\tolr{\tol_{\res}}
\def\tolt{\tol_{\trace}}
\def\toln{\tol_{\nu}}
\def\l{\lambda}
\def\L{\Lambda}
\def\s{\sigma}
\def\S{\Sigma}
\def\R{\mathbb{R}}
\def\C{\mathbb{C}}
\def\Real{\mathrm{Re}}
\def\diag{\mathrm{diag}}
\def\tridiag{\mathrm{tridiag}}

\def\vec{\mathrm{vec}}

\def\trace{\mathrm{trace}}
\def\Span{\mathrm{span}}

\def\twobytwo#1#2#3#4{\bigl[{\hfil#1\atop\hfil#3}{\hfil#2\atop\hfil#4}\bigr]}
\def\twobyone#1#2{\bigl[{\hfil#1\atop\hfil#2}\bigr]}
\def\norm#1{\|#1\|}
\def\normt#1{\|#1\|_2}

\def\ph{\widehat{p}}

\def\mP{\mathcal{P}}
\def\mU{\mathcal{U}}
\def\mM{\mathcal{M}}
\def\mN{\mathcal{N}}
\def\Comp{C}
\def\Comph{A}

\def\Qt{\widetilde{Q}}

\def\DY{\Delta Y}
\def\DW{\Delta W}

\renewcommand{\algorithmicrequire}{\textbf{Input:}}
\renewcommand{\algorithmicensure}{\textbf{Output:}}

\newcommand{\be}{\begin{equation}}
\newcommand{\ee}{\end{equation}}

\begin{document}

\DOI{}
\copyrightyear{}
\vol{}
\pubyear{}
\access{Advance Access Publication Date: Day Month Year}
\appnotes{Paper}
\firstpage{1}


\title[Optimal Damping]
{
Fast Algorithms for Optimal Damping in Mechanical
Systems}

\author{Qingna Li%
        \address{\orgdiv{School of Mathematics and Statistics},
        \orgname{Beijing Key Laboratory on MCAACI,
         and Key Laboratory of Mathematical Theory and Computation in Information
           Security, Beijing Institute of Technology},
        \orgaddress{Beijing, \country{China}}}}
\author{Fran\c{c}oise Tisseur%
       \ORCID{0000-0002-1011-2570}
       \address{\orgdiv{Department of Mathematics},
       \orgname{The University of Manchester},
       \orgaddress{\street{Oxford Road, Manchester}, \postcode{M13 9PL},
       \country{England}}}}


\corresp[*]{Corresponding author:
\href{email:francoise.tisseur@manchester.ac.uk}{francoise.tisseur@manchester.ac.uk}}

\received{8 January 2026}{0}{Year}


\abstract{
Optimal damping consists of determining a vector of damping coefficients $\nu$
that maximizes the decay rate of a mechanical system's response. This problem
can be formulated as the minimization of the trace of the solution of a Lyapunov
equation whose coefficient matrix, representing the system dynamics, depends on
$\nu$. For physical relevance, the damping coefficients $\nu$ must be nonnegative, and the resulting system must be asymptotically stable.
We identify conditions under which the system is never stable or may lose
stability for certain choices of $\nu$. In the latter case, we propose replacing the
constraint $\nu\ge 0$ with $\nu\ge d$, where $d$ is a nonzero nonnegative vector
chosen to guarantee stability.
We derive an
expression for the gradient and Hessian
of the objective function and show that the Karush-Kuhn-Tucker conditions are
equivalent to the vanishing of a nonlinear residual function of $\nu$
at an optimal solution $\nu_*$.
To compute $\nu_*$, we propose a Barzilai-Borwein residual minimization algorithm
(BBRMA), which offers a good balance between simplicity and computational
efficiency but is not guaranteed to converge globally.
We therefore also propose
a spectral projected gradient (SPG) method, which is globally convergent.
The efficiency of both algorithms relies on a fast computation of the gradient
for BBRMA, and both the objective function and its gradient for SPG.
By exploiting the structure of the problem,
we show how to efficiently compute the objective function and its gradient,
with eigenvalue decompositions being the dominant cost in terms of execution time.
Numerical experiments  demonstrate that both BBRMA and SPG require fewer
eigenvalue decompositions than the fast optimal damping algorithm (FODA)
proposed by Jakovčević Stor et al, [\emph{Mathematics}, 10(5):790, 2022.],
and that, although SPG needs extra eigenvalue decompositions when line search is
required, it tends to converge faster than BBRMA, resulting in an overall
lower number of eigendecompositions.}
 \maketitle

\section{Introduction}
\begin{figure}[t]
\begin{center}
\includegraphics[width=.6\textwidth]{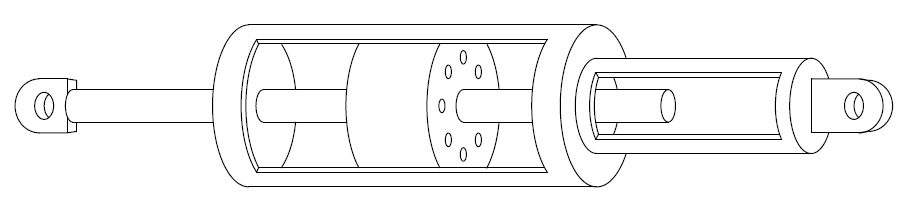}
\caption{Schematic of a viscous damper taken from~\cite[Fig.1]{tasl15}.
It depicts a piston within a cylinder containing a viscous fluid.
As the piston moves, it forces the fluid through narrow passages, such as
orifices, generating a damping force that dissipates energy and reduces
vibrations.}
\label{fig.damper}
\end{center}
\end{figure}

We consider freely oscillating damped vibrational systems
of the form
\be\label{eq-1}
M\ddot{q}(t)+D(\nu)\dot{q}(t)+Kq(t)=0, \quad
q(0) = q_0, \quad \dot{q}(0) = q_1,
\ee
where the $\nbyn$ mass matrix $M$ and the stiffness matrix $K$ are symmetric positive definite,
and the damping matrix $D(\nu)$ is symmetric positive semidefinite.
We assume that $D(\nu)$ can be decomposed as
$$
        D(\nu) = \Dint+ \Dext(\nu),
$$
where the internal damping $\Dint$ is proportionally damped
(i.e., $M,\Dint$ and $K$ are simultaneously diagonalizable
via congruence transformation), and
the external damping $\Dext(\nu)$ is of the form
\begin{equation}\label{def.Dext}
        \Dext(\nu) = \sum_{i = 1}^{k}\nu_i^{}\; D_i D_i^T,
\end{equation}
where
the $D_i\in\R^{n\times r_i}$ are appropriately scaled full rank matrices
that describe the positions of the
$$
        k_d:=\sum_{i=1}^k r_i
$$
dampers in the system, $k$ being the maximum number of distinct linear
viscous dampers,
and the $\nu_i$'s being the damping coefficient of damper $i$.
The above formulation allows for two or more dampers to have the same
damping coefficient,
in which case $r_i\ge 2$ for some $i$. This is particularly useful when dealing
with structures with symmetries. In practice, the total number of dampers $k_d$
is small compared to $n$, i.e., the number of degrees of freedom of the system.
For example, the free oscillations of the London Millennium Footbridge
(see~\cite[Sec .~1]{time01}) is an example of a problem of the form~\eqref{eq-1}.
The bridge
has $k_d=37$
viscous fluid dampers (see Figure~\ref{fig.damper} for a schematic of a viscous
damper) to dissipate the energy of the system \cite{wiki:milleniumbridge}. We
refer to \cite{bhms13} and \cite{vese11} and references therein for other
examples of applications, where such systems arise.

Letting $q(t) = \mathrm{e}^{\l t} x$ in \eqref{eq-1} leads to the quadratic
eigenvalue problem
\begin{equation}\label{def.qep}
Q_\nu(\l)x:=(\l^ 2 M + \l D(\nu) + K)x = 0, \quad x\ne 0,
\end{equation}
where $\l$ is an eigenvalue of $Q_\nu(\l)$ and $x$ is the corresponding eigenvector.
When all the eigenvalues of $Q_\nu(\l)$ strictly lie in the left-hand side of
the complex plane,
the system \eqref{eq-1} is \emph{(asymptotically) stable},
that is, its solution $q(t) = \begin{bmatrix}I&0\end{bmatrix}p(t)$, where
\begin{equation}\label{eq.solq}
         p(t)=\mathrm{e}^{\Comp(\nu)t}p_0,
         \quad p_0=
         \begin{bmatrix}q_0\\ q_1\end{bmatrix},
         \quad
\Comp(\nu) = \begin{bmatrix} 0&I\\ -M^{-1}K & -M^{-1} D(\nu)\end{bmatrix}
\end{equation}
decreases exponentially to zero as $t$ tends to infinity.
The matrix $\Comp(\nu)$ in \eqref{eq.solq} is the companion form of
$M^{-1}Q_\nu(\l)$ and both share the same eigenvalues since
$\Comp(\nu)$ is equivalent to $I_n\oplus Q_\nu(\l)$~\cite{glr09}.
The energy of $q(t)$ at time $t$ is given by
\begin{equation}\label{def.energy}
\mathcal{E}(t) = \frac{1}{2}\dot{q}(t)^T M \dot{q}(t) +
                     \frac{1}{2} q(t)^T K q(t).
\end{equation}
Optimal damping consists of determining a vector $\nu_*\in\R^k_+$ of damping
coefficients that maximizes the rate of decay of $\mathcal{E}(t)$ as $t$ tends to infinity.
This corresponds to the dampers' coefficients that make $q(t)$ converge to zero the fastest
\cite{frla99}, \cite{vese11}.
Different optimization criteria have been considered in the literature.
\begin{itemize}
\item
One criterion is based on the spectral abscissa of $Q_\nu(\l)$ defined by
$\alpha(Q_\nu) = \max \Real (\Lambda(Q_\nu))$,
where $\Lambda(Q_\nu)$ denotes the set of eigenvalues of $Q_\nu(\l)$.
Then for any $\epsilon>0$ there exists a constant $\kappa_\epsilon>0$ such that
$$
        \mathcal{E}(t) \le \kappa_\epsilon  \mathcal{E}(0)
        e^{(2\alpha(Q_\nu)+\epsilon)t}
$$
(see \cite{frla99}, \cite{vese11}).
So we can minimize $\alpha(Q_\nu)$ over all possible values of $\nu$
to maximize the asymptotic rate at which the total energy of the
system decays to zero.
However, this does not control the rate of decay at all times.
Moreover,
$\alpha(Q_\nu)$ can be difficult to compute since it is not a
smooth function of the entries of $Q_\nu(\l)$.

\item Another criterion proposed in~\cite{naki02} and
\cite[Chap.21]{vese11}
consists of finding a $\nu_*\ge 0$ such that
\begin{equation}\label{def.crit}
        \mathcal{E}_T(\nu_*) = \min_{\nu\in\R^k_+} \mathcal{E}_T(\nu),
\end{equation}
where
\begin{equation}\label{eq.ET}
\mathcal{E}_T(\nu) := \int_{\norm{p_0} = 1}
                     \int_0^\infty \mathcal{E}(t) \;dt\; d\s,
\end{equation}
denotes the total energy of the free oscillating system, averaged
over the whole time history
and over all the initial values
$p_0 = \twobyone{q_0}{q_1}$ of unit norm
(to eliminate the dependence on the initial data and time),
$\s$ being a chosen nonnegative measure on the unit sphere in $\R^{2n}$.
\item Other energy-based criteria used by the engineering community are
listed in the survey~\cite{drt19}.
Unlike the criterion~\eqref{def.crit}, they depend on the
external force(s) applied to the system.
\end{itemize}
In this paper, we concentrate on criterion~\eqref{def.crit}.
When $\Comp(\nu)$ is stable or equivalently, \eqref{eq-1} is stable,
it follows from the definition of $\mathcal{E}(t)$ in \eqref{def.energy},
that
$$
\mathcal{E}_T(\nu) = \frac{1}{2}\int_{\norm{p_0} = 1} p_0^T X(\nu) p_0 d\s,
\quad\mbox{where}\quad
X(\nu) :=  \int_0^\infty  \mathrm{e}^{\Comp(\nu)^T t}
      \begin{bmatrix}K&0\\ 0& M\end{bmatrix} \mathrm{e}^{\Comp(\nu)t}\; dt
      \in\R^{2n\times 2n}
$$
is the unique symmetric solution of the Lyapunov equation
\begin{equation}\label{eq.lyap0}
  \Comp(\nu)^T X(\nu) + X(\nu) \Comp(\nu) =
  - \begin{bmatrix}K&0\\ 0& M\end{bmatrix}.
\end{equation}
It is shown in~\cite[Sec.~2]{truh04} that
when $d\s$ is the Lebesgue measure on
the unit sphere $\{p_0\in\R^{2n}: \norm{p_0}~=~1\}$,
$
        \int_{\norm{p_0}=1} p_0^T X(\nu) p_0 d\s
        = \frac{1}{2n}\trace\big(X(\nu)\big)
$.
So the minimization problem \eqref{def.crit} is equivalent to
\begin{equation}\label{def.critb}
        \min_{\nu\in\R^k_+}
        \big\{\trace\big(X(\nu)\big):\;
        \hbox{ $X(\nu)$ solves \eqref{eq.lyap0}},
        \;\
        \Comp(\nu) \hbox{ is stable}
        \big\}.
\end{equation}
Several algorithms have been developed to solve \eqref{def.critb}
when the system is known to be stable for all $\nu$ (see for example \cite{btt11, btt13, jst22, ktt12, ntt13, ttp17, naki02}, and
references therein) but none of these algorithms takes into account the nonnegative constraint
on $\nu$.
As pointed out by Cox et al.~in \cite{cnrv04}, a negative damping coefficient is not
physical as it adds to the energy of the system.
So when this happens, Cox et al. suggest throwing away one damper without
changing the value of the other damping coefficients. This however may not yield an optimal
solution. Moreover, none of these algorithms guarantees that the computed solution
is a strict local minimum.

Our contributions are as follows.
We identify in Section~\ref{sec.prem} conditions under which the system is never stable or may not be stable for
certain values of $\nu$.
In the latter case, we propose to relax the nonnegativity constraint on $\nu$
to $\nu\ge d$ for some nonzero nonnegative vector $d$ chosen such that the system is always stable.

From an optimization perspective, problem \eqref{def.critb} is a nonlinear
optimization problem subject to a lower bound on the vector of damping coefficients. Assuming
the system is stable, we derive expressions in Section~\ref{sec-property} for
the gradient and  Hessian of the objective function. We further show that the
Karush-Kuhn-Tucker (KKT) conditions are equivalent to a nonlinear residual
function dependent on $\nu$ that is equal to zero at an optimum
for~\eqref{def.critb}. To minimize this residual function, we propose in
Section~\ref{sec-algorithm} a Barzilai-Borwein residual minimization algorithm
(BBRMA), which offers a favorable balance between computational efficiency and
simplicity but does not guarantee global convergence. We therefore also
introduce a spectral projected gradient (SPG) method, which, although
potentially more computationally expensive due to the use of a nonmonotone line
search, is globally convergent.

In this work, we focus on small to medium scale problems, including those
arising from model order reduction, where repeated eigenvalue decompositions
remain computationally feasible. For such problems, the practical performance of
the proposed methods critically depends on the efficient evaluation of the
objective function and its derivatives. By exploiting the underlying structure
of the Lyapunov equation and the system matrices, we develop implementations
that significantly reduce the computational cost of these evaluations, with the
dominant expense in terms of execution time being the eigenvalue decompositions.
We show how to inexpensively verify that the computed solution is a strict local
minimizer.
We assess the performance of BBRMA and SPG relative to the fast optimal damping
algorithm (FODA) of Jakovčević Stor et al \cite{jst22}, focusing on the number
of eigenvalue decompositions required to reach convergence. The numerical
experiments reported in Section~\ref{sec.numexp} show that both BBRMA and SPG
require fewer eigenvalue decompositions than FODA while producing nonnegative
damping coefficients. BBRMA converges for the majority of test problems, and
the use of a line search in SPG often leads to faster convergence and a reduced
number of eigenvalue decompositions compared to BBRMA.

\section{Preliminaries}\label{sec.prem}

\subsection{Stability of the system \eqref{eq-1}}\label{sec.stable}

The eigenvalues of the symmetric quadratic eigenvalue problem $Q_\nu(\l)x=0$ in
\eqref{def.qep} with $M$ and $K$ positive definite and
$D(\nu)=\Dint+\Dext(\nu)$ positive semidefinite have nonpositive real
parts~\cite[Thm.~7.1]{lanc66}, \cite{time01}.
Indeed, if $(\l_0,x_0)$ is an eigenpair of $Q_\nu(\l)$. then $\l_0$
is one of the two roots
$(-d_{x_0}\pm\sqrt{d_{x_0}^2-4m_{x_0}k_{x_0}})/{(2m_{x_0})}$ of
the scalar quadratic polynomial
$x_0^HQ_\nu(\l)x_0 = \l^2 m_{x_0} + \l d_{x_0}+k_{x_0}$, where
$m_{x_0}={x_0}^HMx_0^{}>0$, $d_{x_0}=x_0^HD(\nu)x_0^{}\ge 0$,
$k_{x_0}=x_0^HKx_0^{}>0$, and $x_0^H$ denotes the conjugate transpose of $x_0$.
Now,
\begin{itemize}
\item if $\l_0$ is real, then $d_{x_0}^2> 4m_{x_0}k_{x_0}$,
since the nonsingularity of  $K$  implies
$\l_0\ne 0$. So $\l_0<0$;
\item if $\l_0$ is nonreal, then $\Real(\l_0)=-d_{x_0}/(2m_{x_0})\le 0$.
\end{itemize}
Hence the system \eqref{eq-1} is not stable if and only if $d_{x_0}=0$ for some
eigenvector
$x_0$ of $Q_\nu(\l)$.
So a sufficient condition for stability is for
$\Dint$ to be nonsingular.
Fortunately, in many applications, $\Dint$ corresponds to the Rayleigh damping
matrix
\begin{equation}\label{eq-cint2}
        \Dint = \alpha M + \beta K,
        \quad\alpha,\beta\ge 0,
\end{equation}
with $\alpha,\beta$ not both zero, or the critical damping matrix \cite{inma01}
\begin{equation}\label{eq-cint1}
\Dint = \alpha M^{\frac12}\sqrt{M^{-\frac12}KM^{-\frac12}}M^{\frac12},
\quad \alpha>0,
\end{equation}
 both of which are positive definite so that $d_{x_0}\ne 0$ and
the system \eqref{eq-1} is stable for all $\nu$.

We need further results when there is no internal damping or $\Dint$ is singular.
\begin{theorem}\label{thm.stable}
Let $\nu\in\R_+^k$, $\nu\ne 0$.
The system~\eqref{eq-1} is stable if and only if
$(x^H\Dint x,x^H\Dext(\nu)x)\ne (0,0)$
for any eigenvector $x$ of $Q_\nu(\l)$ in~\eqref{def.qep}.
\end{theorem}
\vspace*{-.5cm}
\begin{proof}
It is easy to see that the system is stable if and only if
$x^HD(\nu)x>0$ for all the eigenvectors $x$ of $Q_\nu(\l)$.
The result follows
since both $\Dint$ and $\Dext(\nu)$ are positive semidefinite.
\end{proof}
As mentioned in the introduction, $\Dint$ is such that
$M,\Dint,K$ are simultaneously diagonalizable, or equivalently, the
Caughey-O'Kelly commutativity condition $\Dint M^{-1}K=KM^{-1}\Dint$
holds~\cite{caok65}.
Let $\Phi$ denotes an $\nbyn$ nonsingular matrix that
simultaneously diagonalizes $M$ and $K$,
\begin{equation}\label{def.PhiOm}
\Phi^T M \Phi =I_n,
\quad
\Phi^T K \Phi = \Omega^2,
\quad
\Omega = \diag(\omega_1,\dots, \omega_n),
\quad
0< \omega_1 \le\dots\le \omega_n.
\end{equation}
Then $\Dint M^{-1}K=KM^{-1}\Dint$ implies that
$ \Phi^T \Dint \Phi$ commutes with $\Omega^2$.
Let $\Omega^2 = \diag(\kappa_1 I_{\ell_1},\ldots,\kappa_r I_{\ell_r})$, where
the $\kappa_i$ are distinct with $\kappa_i = \omega_j^2$ for some $j$ and
$\ell_i$ is the multiplicity of $\kappa_i$.
Then $\Phi^T \Dint \Phi$ is block diagonal with diagonal blocks of size
$\ell_1, \ldots,\ell_r$ and since these blocks are symmetric, there exists
a block diagonal matrix $\Phi_1$ whose $i$th block is orthogonal
and diagonalizes the $i$th block of $\Phi^T \Dint \Phi$ and leaves
$\Phi^T M\Phi$ and $\Phi^T K\Phi$ unchanged.
In what follows we let $\Phi\equiv\Phi\Phi_1$ and
\begin{equation}\label{def.Gamma}
\Gamma := \Phi^T \Dint \Phi=\diag(\gamma_1,\ldots,\gamma_n).
\end{equation}
If $\Dint$ is singular, then by Sylvester's law of inertia,
some of the $\gamma_i$'s are equal to zero.
Let
\begin{equation}\label{def.Sigmanu}
R := \Phi^T\begin{bmatrix} D_1 & \dots & D_k\end{bmatrix}\in\R^{n\times k_d},
\quad
\Sigma_\nu := \diag(\nu_1 I_{r_1},\dots,\nu_k I_{r_k})\in\R^{k_d\times k_d},
\quad
R_\nu := R\Sigma_\nu^{1/2}
\end{equation}
so that $\Phi^T\Dext(\nu)\Phi^{} = R_\nu^{}R_\nu^T$.
Define
$$
         [n]:=\{1,\dots,n\}.
$$
We have that:
\begin{enumerate}
\item[(i)]
If $\gamma_j=0$ and $e_j^TR=0$  for some $j\in[n]$,
then the system \eqref{eq-1} is not stable for all $\nu$.
\item[(ii)]
If $\gamma_j=0$ and $e_j^TR_\nu =0$ for some $j\in[n]$ and
some $\nu\in\R_+^k$, then the system~\eqref{eq-1} is not stable for this particular $\nu$.
\end{enumerate}
Case (ii) can only happen when some of the damping coefficients are equal to zero.
We then suggest replacing the nonnegative constraint on $\nu$
with $\nu\ge d$ for some nonzero nonnegative vector $d$ chosen such that
for any such $\nu$, $e_j^TR_\nu \ne 0$. We will discuss later how to deal with
this type of constraint.

\subsection{Restating the optimization problem in the modal
basis}\label{sec.modal}


In practical applications, one is usually
interested in damping only a certain part of the
spectrum, say, all  eigenvalues of the undamped quadratic $\l^2 M+K$ that are
smaller than $\omega_{\max}$ or the eigenmodes that affects most the
response of the system under the influence of some forces, such as
seismic forces \cite{ztk19}.
To this end, it is convenient to rephrase the problem in the modal basis
spanned by the columns of the $\nbyn$ nonsingular matrix $\Phi$
in~\eqref{def.PhiOm}, and to order the
columns of $\Phi$ and the diagonal entries of
$\Omega$ such that
$(\omega_i,\phi_i)$, $i\in[s]$ correspond to the $s$ modes to dampen.

Rewrite the energy $\mathcal{E}(t)$ in \eqref{def.energy} as
$\mathcal{E}(t) = \frac{1}{2} \ph(t)^T\ph(t)$, where
$\ph(t) = \mathrm{e}^{\Comph(\nu)t}\ph_0$ with $\ph_0 = T_\Phi^{-1}p_0$,
and
\begin{equation}\label{def.Comph}
        T_\Phi =\begin{bmatrix}\Phi\Omega^{-1} &0\\ 0&\Phi\end{bmatrix},
        \qquad
        \Comph(\nu)= T_\Phi^{-1}\Comp(\nu) T_\Phi=
        \begin{bmatrix} 0&\Omega\\ -\Omega & -\Phi^TD(\nu)\Phi\end{bmatrix}.
\end{equation}
Then $\mathcal{E}_T(\nu)$ in \eqref{eq.ET} becomes
$$
\mathcal{E}_T(\nu) = \frac{1}{2}\int_{\norm{\ph_0} = 1} \ph_0^T W(\nu) \ph_0 d\s,
$$
where
$W(\nu) =  \int_0^\infty \mathrm{e}^{\Comph(\nu)^Tt}
             \mathrm{e}^{\Comph(\nu)t} dt
$
is the unique
symmetric positive definite solution to the Lyapunov equation
\begin{equation}\label{eq.lyap1}
        \Comph(\nu)^T W(\nu) + W(\nu) \Comph(\nu) = - I_{2n}.
\end{equation}
Let $\s=\s_1\times\s_2\times\s_1\times \s_2$, where
$\s_1$ is a measure on $\Span\{\phi_1,\dots,\phi_s\}$ generated by a Lebesgue
measure on $\R^{2n}$ and $\s_2$ is a Dirac measure on
$\Span\{\phi_{s+1},\dots,\phi_n\}$.
It is shown in \cite{truh04} that
\begin{equation}\label{def.G}
        \frac{1}{2}\int_{\norm{\ph_0} = 1} \ph_0^T W(\nu) \ph_0 d\s
        = \trace\big(W(\nu)Z\big),
        \qquad
        Z = \frac{1}{2s}
          GG^T,
          \qquad
          G=
         \begin{bmatrix}I_s&0\\0&0\\ 0&I_s\\0&0\end{bmatrix}\in \R^{2n\times
         2s}.
\end{equation}
Since $\trace\big(W(\nu)Z\big)=\frac{1}{2s}\trace\big(G^T W(\nu) G\big)$, the
minimization problem in \eqref{def.critb} is equivalent to
\begin{equation}\label{def.critbb}
        \min_{\nu\in\R^k_+}
        \big\{\trace(G^TW(\nu) G):
        \ \mbox{$W(\nu)$ solves~\eqref{eq.lyap1}},
        \ \mbox{$\Comph(\nu)$ is stable}
        \big\}.
\end{equation}
Also,
\begin{equation}\label{eq.traceWY}
\trace\big(W(\nu) Z\big) = \trace \Big(\int_0^\infty \mathrm{e}^{\Comph(\nu)^Tt}
             \mathrm{e}^{\Comph(\nu)t} dt Z\Big)
             =\trace \Big(\int_0^\infty \mathrm{e}^{\Comph(\nu)t}Z
             \mathrm{e}^{\Comph(\nu)^Tt} dt\Big)
             =\trace\big(Y(\nu)\big),
\end{equation}
where $Y(\nu)$ solves the Lyapunov equation
\begin{equation}\label{eq.lyap2}
        \Comph(\nu) Y(\nu) + Y(\nu) \Comph(\nu)^T = - Z.
\end{equation}
So \eqref{def.critbb} is equivalent to
\begin{equation}\label{def.critbb_dual}
        \min_{\nu\in\R^k_+}
        \big\{\trace\big(Y(\nu) \big):
        \ \mbox{$Y(\nu)$ solves~\eqref{eq.lyap2}},
        \ \mbox{$\Comph(\nu)$ is stable}
        \big\}.
\end{equation}

\subsection{A Toy Example}\label{sec.toyexample}
We consider the small problem in~\cite[Example~6]{cnrv04} defined by
$$
M = I_2,
\quad \Dint = 0,
\quad
\Dext(\nu)= \nu_1 e_1^{}e_1^T +\nu_2(e_2^{}-e_1^{})(e_2^{}-e_1^{})^T,
\quad
K = \begin{bmatrix}1 &-1\\ -1 & 201\end{bmatrix}.
$$
The corresponding system is stable for all $\nu>0$.
Now if $\nu_1=0$ and $\nu_2\ne 0$, then $x^T \Dext x=0$ if and only if $x$ is a multiple of
$\twobyone{1}{1}$. But the latter is not an eigenvector of $Q_\nu(\l)=\l^2 I + \l \Dext+K$.
So by Theorem~\ref{thm.stable}, the system is stable.
Similarly, if $\nu_1\ne 0$ and $\nu_2=0$, then $x^T \Dext x=0$ if and only if
$x$ is a multiple of $e_2$. But the latter is not an eigenvector of $Q_\nu(\l)$, so
the system is stable.

Solving~\eqref{def.critb} or equivalently \eqref{def.critbb_dual}
with $Z=I/4$ and no nonnegative constraints on $\nu$ leads to $\nu_*\approx\twobyone{-2.59}{4.75}$
with $\trace\big(Y(\nu_*)\big) \approx 0.67$.
So, the first damper has negative damping coefficient and Cox et al.~\cite{cnrv04} suggest setting it to zero.
Now for $\widehat{\nu}_*=\twobyone{0}{4.75}$, we find that
$\trace\big(Y(\widehat\nu_*)\big) \approx 1.19$.
But $\nu_{**}=\twobyone{0}{2.72}\ge 0$ is a better minimizer
since  $\trace(Y(\nu_{**})) \approx 0.73<\trace\big(Y(\widehat\nu_*)\big)$.
In what follows, we show how to compute such a solution.

\section{Gradient, Hessian, and Optimality Conditions}\label{sec-property}

Before stating the
KKT conditions
associated with the minimization problem~\eqref{def.critbb_dual},
we derive the gradient and Hessian of the objective function
\begin{equation}\label{def.f}
f(\nu)=\trace(Y(\nu)),
\end{equation}
where $Y(\nu)$ solves~\eqref{eq.lyap2}.
We need the following results from
\cite[p.491, Eq.(6.5.7) and Eq.(6.5.8)]{hojo91}.
\begin{lemma}\label{lem.diffA}
If the $\nbyn$ matrix $F(x)$ is differentiable and nonsingular at some
point $x\in\R$, then
\begin{align}
\frac{d}{dx}F(x)^{-1} &= -F(x)^{-1}\left(\frac{d}{dx}F(x)\right) F(x)^{-1}.
\label{eq.dAdt}        \\
\frac{d}{dx}\trace\big(F(x)\big) &=
\trace\left( \frac{d}{dx}F(x)\right). \label{eq.dtracedt}
\end{align}
\end{lemma}
We interchangeably use the matrix equality constraint \eqref{eq.lyap2}
or its vector form
\begin{equation}\label{def.mP}
\bigl(
\underbrace{I_{2n}\otimes \Comph(\nu) + \Comph(\nu)\otimes  I_{2n}}_{\mP(\nu)}\bigr)
\vec(Y(\nu))=-\vec(Z)
\end{equation}
obtained by applying the $\vec$ operator, which stacks the columns of a matrix into
one long vector. Here $\otimes$ denotes the Kronecker product between two
matrices.
The next result follows from the fact that for
two eigenpairs $(\l_i,x_i)$, $i=1,2$, of $\Comph(\nu)$,
$(\l_1+\l_2, x_1\otimes x_2)$ is an eigenpair of $\mP(\nu)$.
\begin{lemma}\label{lem.Pnonsing}
Let $\nu\in\R^k_+$. Then $\Comph(\nu)$ is stable for some $\nu\in\R^k_+$
if and only if $\mP(\nu)=I_{2n}\otimes \Comph(\nu) + \Comph(\nu)\otimes  I_{2n}$ is
nonsingular.
\end{lemma}
Let us rewrite $\Comph(\nu)$ in~\eqref{def.Comph} as
\begin{equation}\label{def.Anu}
\Comph(\nu) = \Comph(0)
-
\sum_{i=1}^k \nu_i U_i^{}U_i^T,
\quad
U_i =
\begin{bmatrix} 0\\ R_i\end{bmatrix}\in\R^{2n\times r_i},
\quad
R_i = \Phi^T D_i,
\quad i\in[k]
\end{equation}
so that
\begin{equation}\label{eq.dAdnu}
\frac{\partial}{\partial \nu_i} \mP(\nu) =  -(I\otimes U_i^{}U_i^T +
U_i^{}U_i^T\otimes I), \quad i\in[k].
\end{equation}
\begin{theorem}\label{thm.nablaf}
Let $\nu\in\R^k$ be such that $\Comph(\nu)$ is stable,
$f$ be as in \eqref{def.f}, and
$W(\nu)$ be the solution to~\eqref{eq.lyap1}.
Then
\begin{align}
\nabla f(\nu) &= -\mU(Y(\nu),W(\nu)), \label{eq.nablaf}\\
\nabla^2 f(\nu) &= -\left(\mM(\nu) \mP(\nu)^{-T} \mN(\nu)^T
                 +\mN(\nu) \mP(\nu)^{-1} \mM(\nu)^T\right),
                 \label{eq.nabla2f}
\end{align}
with $\mathcal{U}: \mathbb R^{2n\times 2n}\times \mathbb R^{2n\times 2n} \to\R^{k}$
such that
\begin{equation}\label{def.mU}
e_i^T\mathcal{U}(Y(\nu),W(\nu)) =
\trace\big(W(\nu)^T(U_i^{}U_i^T Y(\nu) + Y(\nu)U_i^{}U_i^T)\big),
\quad i\in [k],
\end{equation}
$U_i$ as in~\eqref{def.Anu}, $\mP$ as in~\eqref{def.mP},
and
$\mathcal{M},\mathcal{N}\in\R^{k\times 4n^2}$ such that
$$
     e_i^T\mathcal{M}(\nu) = \vec(U_i^{}U_i^TY(\nu) + Y(\nu) U_i^{} U_i^T)^T,
     \quad
     e_i^T\mathcal{N}(\nu) = \vec(W(\nu) U_i^{}U_i^T + U_i^{}U_i^TW(\nu))^T,
     \quad i\in[k].
$$
\end{theorem}
\begin{proof}
On using~\eqref{eq.dtracedt} and the fact that
for two matrices $A$ and $B$ of appropriate sizes,
\begin{equation}\label{identity-tracevec}
\trace(A^TB) = \vec(A)^T\vec(B),
\end{equation}
we find that for $i\in[k]$,
$$
\frac{\partial }{\partial \nu_i}f(\nu)
=
\frac{\partial }{\partial \nu_i} \trace\big(Y(\nu)\big)
=\trace\left(\frac{\partial }{\partial \nu_i}Y(\nu)\right)
=  \vec(I)^T\vec\left(\frac{\partial }{\partial \nu_i}Y(\nu)\right)
=\vec(I)^T\frac{\partial }{\partial \nu_i}\vec\left(Y(\nu)\right).
$$
It follows from \eqref{def.mP} and Lemma~\ref{lem.Pnonsing} that
$\vec\big(Y(\nu)\big) = -\mP(\nu)^{-1} \vec(Z)$, and
on using~\eqref{eq.dAdt} and \eqref{eq.dAdnu}, we find that
\begin{equation}\label{eq.dmPdnui}
\frac{\partial }{\partial \nu_i}\mP(\nu)^{-1} =
\mP(\nu)^{-1}(I\otimes U_i^{}U_i^T +  U_i^{}U_i^T\otimes I)\mP(\nu)^{-1}
\end{equation}
so that
\begin{equation}\label{eq.dfdnui}
\frac{\partial }{\partial \nu_i}f(\nu)
= -
\vec(I)^T\mP(\nu)^{-1}
(I\otimes U_i^{}U_i^T +  U_i^{}U_i^T\otimes I)
\mP(\nu)^{-1}\vec(Z).
\end{equation}
Now, $\vec\big(W(\nu)\big) = -\mP(\nu)^{-T} \vec(I)$ so that
\begin{align*}
\frac{\partial }{\partial \nu_i}f(\nu)
&= -
\vec\big(W(\nu)\big)^T
(I\otimes U_i^{}U_i^T +  U_i^{}U_i^T\otimes I)
\vec\big(Y(\nu)\big)\\
&= -
\vec\big(W(\nu)\big)^T
\vec\big(U_i^{}U_i^T Y(\nu)+ Y(\nu) U_i^{}U_i^T\big) \\
&= -\trace\big(W(\nu)^T(U_i^{}U_i^T Y(\nu)+ Y(\nu) U_i^{}U_i^T)\big),
\end{align*}
where we used~\eqref{identity-tracevec} for the last equality.
The formula for $\nabla f(\nu)$ in \eqref{eq.nablaf} follows.

On using~\eqref{eq.dmPdnui} and the product rule, we find that
\begin{align*}
\frac{\partial}{\partial \nu_i}
\mP(\nu)^{-1}
(I\otimes U_j^{}U_j^T +  U_j^{}U_j^T\otimes I)
\mP(\nu)^{-1}
&=
\mP(\nu)^{-1}
(I\otimes U_i^{}U_i^T +  U_i^{}U_i^T\otimes I)
\mP(\nu)^{-1}\times\\
&
(I\otimes U_j^{}U_j^T +  U_j^{}U_j^T\otimes I)
\mP(\nu)^{-1}\\
&+
\mP(\nu)^{-1}
(I\otimes U_j^{}U_j^T +  U_j^{}U_j^T\otimes I)\times\\
&\mP(\nu)^{-1}
(I\otimes U_i^{}U_i^T +  U_i^{}U_i^T\otimes I)
\mP(\nu)^{-1}.
\end{align*}
Hence, starting with~\eqref{eq.dfdnui} and using the above as well as
$\vec\big(W(\nu)\big) = -\mP(\nu)^{-T} \vec(I)$ and
$\vec\big(Y(\nu)\big) = -\mP(\nu)^{-1} \vec(Z)$, we have that
\begin{align*}
\frac{\partial^2 }{\partial \nu_i\partial \nu_j}f(\nu)
&=- \vec\big(W(\nu)\big)^T
(I\otimes U_i^{}U_i^T +  U_i^{}U_i^T\otimes I)
\mP(\nu)^{-1}
(I\otimes U_j^{}U_j^T +  U_j^{}U_j^T\otimes I)
\vec\big(Y(\nu)\big)
\\
&{}\qquad-
\vec\big(W(\nu)\big)^T
(I\otimes U_j^{}U_j^T +  U_j^{}U_j^T\otimes I)
\mP(\nu)^{-1}
(I\otimes U_i^{}U_i^T +  U_i^{}U_i^T\otimes I)
\vec\big(Y(\nu)\big)\\
&=
-\vec\big(U_i^{}U_i^TW(\nu))+W(\nu)U_i^{}U_i^T\big)^T
\mP(\nu)^{-1}
\vec\big(U_j^{}U_j^T Y(\nu)+ Y(\nu) U_j^{}U_j^T\big)
\\
&{}\qquad
-\vec\big(U_j^{}U_j^TW(\nu)+W(\nu)U_j^{}U_j^T\big)^T
\mP(\nu)^{-1}
\vec\big(U_i^{}U_i^T Y(\nu)+ Y(\nu) U_i^{}U_i^T\big)\\
&=
-e_i^T \mN (\nu)\mP(\nu)^{-1}\mM (\nu)^T e_j
-e_j^T \mN(\nu)\mP(\nu)^{-1}\mM (\nu)^T e_i,
\end{align*}
which completes the proof.
\end{proof}

Let
\begin{equation}\label{def.Lagrange}
\mathcal L(\nu,\mu)= \trace\big(Y(\nu)\big)-\nu^T\mu,
\end{equation}
be the Lagrange function associated with~\eqref{def.critbb_dual}, where
the entries of $\mu\in\R^k$ are the Lagrange multipliers.
The KKT conditions for~\eqref{def.critbb_dual} are given by
\begin{subnumcases}{\label{def.kkt_cond}}
\nabla_\nu \mathcal L(\nu,\mu)= 0,\label{kkt.stationary2}\\
\nu\ge 0, \label{kkt.primal}\\
\mu\ge 0,\ \nu^T\mu=0.
\label{kkt.comp.dual}
\end{subnumcases}
with
the stationary condition in~\eqref{kkt.stationary2},
the primal feasibility condition in~\eqref{kkt.primal},
the dual feasibility and complementary slackness condition
in~\eqref{kkt.comp.dual}.
It follows from~\eqref{def.Lagrange} that~\eqref{kkt.stationary2} is equivalent to
$\nabla f(\nu)=\mu$, which we
use to eliminate $\mu$ from the KKT system~\eqref{def.kkt_cond}.
This leads to the reduced KKT conditions
$$
\nabla f(\nu) \ge 0,\ \nu\ge 0,\ \nu^T \nabla f(\nu) = 0,
$$
which are equivalent to
\begin{equation}\label{def.h}
h(\nu):=\nu-\max(\nu-\nabla f(\nu),0)=0.
\end{equation}
We will refer to $h(\nu)$ as the residual function.

\section{Two Algorithms for Solving Problem \eqref{def.critbb_dual}}\label{sec-algorithm}

The nonlinear system~\eqref{def.h} is nonsmooth and therefore cannot be addressed directly by standard smooth Newton-type methods.
A smoothing approach is considered in~\cite{wlt26}. Other globalized algorithms for nonsmooth systems could also be considered.
Here we consider a different approach.

In their seminal article on two-point step size gradient methods, Barzilai and
Borwein note that the iteration
\begin{equation}\label{def.BBiter_init}
 \nu^{(j+1)} = \nu^{(j)}-\eta_{BB}^{(j)}h^{(j)},
\end{equation}
where  $h^{(j)} = h\big(\nu^{(j)}\big)$
and stepsize $\eta_{BB}^{(j)}$ given by
\begin{equation}\label{def.BBstep}
\eta_{BB}^{(j)} =  \frac{\norm{\nu^{(j)}-\nu^{(j-1)}}^2}{
(\nu^{(j)}-\nu^{(j-1)})^T(h^{(j)}-h^{(j-1)})},
\end{equation}
is applicable to the solution of the nonlinear equation $h(\nu)=0$
(see \cite[Sec.~2]{babo88}).
So we apply~\eqref{def.BBiter_init} to our residual function
$h(\nu)$ in~\eqref{def.h} but to maintain the nonnegativity of $\nu$, we replace \eqref{def.BBiter_init} with
\begin{equation}\label{def.BBiter}
 \nu^{(j+1)} = \max(\nu^{(j)}-\eta_{BB}^{(j)}h^{(j)},0).
\end{equation}
This leads to the BB-step residual minimization algorithm (BBRMA) displayed in
Algorithm~\ref{alg.BBRMA}, where $\epsilon$ is some given tolerance on the norm of the residual function.
\begin{algorithm}[th]
\caption{BB-step Residual Minimization Algorithm (BBRMA)}\label{alg.BBRMA}
\renewcommand{\algorithmicrequire}{\textbf{Input:}}
\renewcommand{\algorithmicensure}{\textbf{Output:}}
\begin{algorithmic}[1]
     \Require $\nu^{(0)}\in\R^k_+$, $\eta_{BB}^{(0)}>0$, $\epsilon>0$,
     gradient function $\nabla f$ in~\eqref{eq.nablaf}
     \State $j=0$,
     $\nabla f^{(0)}=\nabla f(\nu^{(0)})$,
     $h^{(0)} = \nu^{(0)}-\max(\nu^{(0)}-\nabla f^{(0)},0)$
     \While{$\norm{h^{(j)}}>\epsilon$}
         \State If $j>0$ then compute the BB-stepsize $\eta_{BB}^{(j)}$ in
         \eqref{def.BBstep}.

         \State  $\nu^{(j+1)} = \max\big(\nu^{(j)}-\eta_{BB}^{(j)}h^{(j)},0\big)$
              \label{step.eta-nu}
         \State $\nabla f^{(j+1)}=\nabla f(\nu^{(j+1)})$
         \State $h^{(j+1)}=\nu^{(j+1)}-\max(\nu^{(j+1)}
         -\nabla f^{(j+1)},0)$
         \State $j=j+1$
     \EndWhile
     \Ensure $\nu=\nu^{(j)}$
    \end{algorithmic}
\end{algorithm}

The Barzilai–Borwein (BB) method for solving nonlinear equations typically exhibits fast local
convergence, but is not guaranteed to be globally convergent.
Since the feasible set $\R_+^k$ is convex, we may instead apply the globally convergent spectral projected gradient (SPG) algorithm developed in~\cite[Alg.~2.2]{bmr00} to the problem $\min_{\nu\in\R^{k}_+} f(\nu)$. 
The SPG method combines the classical projected gradient approach~\cite{gold64} with a global BB (spectral) nonmonotone line search scheme~\cite{Rayd97}.
With the SPG framework, the vector of damping coefficients is updated according to
$$
\nu^{(j+1)}=\nu^{(j)}+ \alpha^{(j)}d^{(j)},
$$
where step length $\alpha^{(j)}$ is determined by the nonmonotone line search. The search direction is given by
$$
d^{(j)} = \max\big(\nu^{(j)}-\eta^{(j)}\nabla f^{(j)},0\big)-\nu^{(j)},
$$
with $\eta^{(j)}$ denoting the (bounded) BB step length computed to minimize the objective function $f(\nu)$.
This leads to Algorithm~\ref{alg.SPG}.


\begin{algorithm}[h]
\caption{Spectral projected gradient algorithm (SPG)}\label{alg.SPG}
\renewcommand{\algorithmicrequire}{\textbf{Input:}}
\renewcommand{\algorithmicensure}{\textbf{Output:}}
    \begin{algorithmic}[1]
         \Require
          $\nu^{(0)}\in\R^k_+$, $\eta^{(0)}>0$, $\epsilon>0$,
          $\sigma,\rho\in (0,1)$, $0<\eta_{min}<\eta_{max}$, integer $M_0\ge1$,
          function $f$ and its gradient $\nabla f$
         \State $j=0$, $\nabla f^{(0)}=\nabla f(\nu^{(0)})$,
          $h^{(0)} = \nu^{(0)}-\max(\nu^{(0)}-\nabla f^{(0)},0)$
         \While{$\norm{h^{(j)}}>\epsilon$}
           \If {$j>0$}
           \State $p = \big(\nu^{(j)}-\nu^{(j-1)}\big)^T
             \big(\nabla f^{(j)}-\nabla f^{(j-1)}\big)$
           \If{$p\le 0$ }
           \State $\eta^{(j)} =\eta_{\max}$
           \Else{}
           \State
           $\eta_{BB}^{(j)} =\norm{\nu^{(j)}-\nu^{(j-1)}}^2/p$
           \label{step.BB}
           \State $\eta^{(j)} =
           \max\big\{\eta_{\min},\min\{\eta_{BB}^{(j)},\eta_{\max}\}\big\}$
           \EndIf
         \EndIf
         \State $d^{(j)} =
         \max\big(\nu^{(j)}-\eta^{(j)}\nabla f^{(j)},0\big)-\nu^{(j)}$
         \label{step.direction}
        \State Find $\alpha^{(j)} = \rho^{m_j}$, where $m_j$ is the
        smallest nonnegative integer such that
       \[
       f(\nu^{(j)}+\alpha^{(j)}d^{(j)})\le
       f_{\max}+\sigma\alpha^{(j)}d^{(j)T}\nabla f^{(j)},
       \]
       \quad\ \
       with $f_{\max} = \max\big\{f(\nu^{(j-\ell)}) \ :\ \ 0\le \ell\le
         \min\{j,M_0-1\}\big\}$.
       \label{step.linesearch}
        \State  $\nu^{(j+1)}=\nu^{(j)}+ \alpha^{(j)}d^{(j)}$\label{step.update}
        \State  $\nabla f^{(j+1)} = \nabla f(\nu^{(j+1)})$
        \State $h^{(j+1)}
         =\nu^{(j+1)}-\max(\nu^{(j+1)}-\nabla f^{(j+1)},0)$
        \State $j=j+1$
      \EndWhile
      \Ensure $\nu=\nu^{(j)}$.
    \end{algorithmic}
\end{algorithm}

\subsection{Solving the Lyapunov Equations in \eqref{eq.lyap1} and \eqref{eq.lyap2}}

Algorithm~\ref{alg.BBRMA} requires an
evaluation of $\nabla f(\nu)$ at each iteration, whereas
Algorithm~\ref{alg.SPG} requires both $\nabla f(\nu)$ and $f(\nu)$.
Since from Theorem~\ref{thm.nablaf},
\begin{equation}\label{eq.nablaf_i}
e_i^T\nabla f(\nu) =
-\trace\big(W(\nu)^T(U_i^{}U_i^T Y(\nu) + Y(\nu)U_i^{}U_i^T)\big)
= -2\trace\big(U_i^T Y(\nu)W(\nu)^T U_i\big),
\end{equation}
any efficient implementation of Algorithm~\ref{alg.BBRMA} and Algorithm~\ref{alg.SPG}
requires an efficient solution of the Lyapunov equations \eqref{eq.lyap1} and~\eqref{eq.lyap2}.
For this, we need to assume
that $\Comph(\nu)$ is diagonalizable so that there exists $T\in\C^{2n\times
2n}$ nonsingular such that
\begin{equation}\label{eq.diag}
        T^{-1}\Comph(\nu) T^{} = \L=\diag(\l_1,\dots,\l_{2n}).
\end{equation}
Note that $T$ depend on $\nu$ but to avoid more clutter with the notation, we
write $T$ in place of $T(\nu)$ or $T_\nu$.
The same applies to  $\L$ and the $\l_j$.
Since $\Comph(\nu)$ is stable, it follows from~\cite[Sec.~5.0.4]{hojo91} that
\begin{equation}\label{eq.solYW}
        Y(\nu) =
        -\frac{1}{2s} T\left(L(\overline{\L})\circ (T^{-1}G)(T^{-1}G)^H\right)T^H,
        \qquad
        W(\nu)  = -T^{-H}\left(L(\L)\circ T^HT\right)T^{-1},
\end{equation}
where $\overline{\Lambda}$ denotes the conjugate of $\Lambda$,
\begin{equation}\label{eq.L}
\big(L(\L)\big)_{ij} = (\bar{\l_i}+\l_j)^{-1}
\end{equation}
and $\circ$ denotes the Hadamard product.

Although $Y(\nu)$ and $W(\nu)$ in \eqref{eq.solYW} can be kept in
factored form for the evaluation of $f(\nu)$ and $\nabla f(\nu)$,
we nevertheless need to form
the two matrix-matrix products $(T^{-1}G)(T^{-1}G)^H$ and $T^HT$.
When the number $s$ of modes to dampen is small compare to the size $n$ of the
system, we can exploit the low rank property of
$G\in\R^{2n\times 2s}$ to reduce the cost of forming $(T^{-1}G)(T^{-1}G)^H$ but
unfortunately not that of $T^HT$.
So evaluating $\nabla f(\nu)$ at each iteration in this way is expensive.

A way around this is to decompose $Y(\nu)$ and $W(\nu)$ as
\begin{equation}\label{def.W0}
      Y(\nu) := Y_0 + \Delta Y(\nu),
      \qquad
      W(\nu) := W_0 + \Delta W(\nu),
\end{equation}
where $Y_0$ and $W_0$ are the unique symmetric solutions to
the Lyapunov equations~\eqref{eq.lyap2} and~\eqref{eq.lyap1} for $\nu=0$,
respectively,
assuming that $A(0)$ is stable and diagonalizable.
Then
$\Upsilon:=(\Gamma^2-4\Omega^2)^{1/2}$ is nonsingular since in this case,
$\gamma_i=2\omega_i$ for some $i$, which implies that
$-\omega_i$ is a defective eigenvalues of $Q_\nu(\l)$ and hence of $A(0)$,
contradicting the assumption that $A(0)$ is diagonalizable.
We then have an explicit expression for its eigendecomposition,
$A(0)=\twobytwo{0}{\Omega}{-\Omega}{-\Gamma}= T_0^{}\Lambda_0 T_0^{-1}$, that is
given by
$$
\Lambda_0 = \begin{bmatrix} -\frac{1}{2}(\Gamma+\Upsilon)&\\
                        &-\frac{1}{2}(\Gamma-\Upsilon)
    \end{bmatrix},
    \quad
T_0= \begin{bmatrix}
-\frac{1}{2}(\Gamma-\Upsilon)\Omega^{-1}&
-\frac{1}{2}(\Gamma+\Upsilon)\Omega^{-1}
\\
I_n&I_n
\end{bmatrix},
$$
$$
T_0^{-1}= \begin{bmatrix}
\Omega\Upsilon^{-1} & \frac{1}{2}(\Gamma+\Upsilon)\Upsilon^{-1}\\
-\Omega\Upsilon^{-1} & -\frac{1}{2}(\Gamma-\Upsilon)\Upsilon^{-1}\\
\end{bmatrix}.
$$
Since $T_0$ and $T_0^{-1}$ are block $2\times 2$ matrices with $\nbyn$ diagonal
blocks, it follows from \eqref{eq.solYW} with $\nu=0$ that
\begin{equation}\label{eq.Y0W0}
Y_0 = \begin{bmatrix}
\Upsilon_1^{} & \Upsilon_2\\ \Upsilon_2^{} & \Upsilon_3^{}
\end{bmatrix},
\qquad
W_0 =\begin{bmatrix}
\Psi_1^{} & \Psi_2\\ \Psi_2^{} & \Psi_3^{}
\end{bmatrix}
\end{equation}
are also block $2\times 2$  with $\nbyn$ diagonal
blocks $\Upsilon_i$, $\Psi_i$, $i=1,2,3$
that can be computed in $O(n)$ operations as the sum and product of diagonal matrices.
Note that a $\nu$ different than $0$ can be chosen in the expansion
\eqref{def.W0} when  $A(0)$ is not stable and forming
$Y_0$ and $W_0$ will be more costly, but this will only need to be done once.

Letting
\begin{equation}\label{def.Unu}
     U:=\begin{bmatrix}U_1& \ldots &U_k\end{bmatrix}\in\R^{2n\times k_d},
     \quad
     U_\nu := U\S_\nu^{1/2}=\begin{bmatrix}0\\R_\nu\end{bmatrix}\in\R^{2n\times k_d}
\end{equation}
with $U_i$ as in \eqref{def.Anu} and $\S_\nu,R_\nu$ as in~\eqref{def.Sigmanu},
we rewrite~\eqref{eq.lyap2} and~\eqref{eq.lyap1} as
\begin{align}
\Comph(\nu) \DY(\nu) + \DY(\nu) \Comph(\nu)^T &=
U_\nu^{}U_\nu^TY_0+ Y_0U_\nu^{} U_\nu^{T},\label{eq.lyap2bis}
\\
\Comph(\nu)^T \DW(\nu) + \DW(\nu) \Comph(\nu)
&= U_{\nu}^{}U_\nu^T W_0+ W_0 U_\nu^{} U_{\nu}^{T}.
\label{eq.lyap1bis}
\end{align}
So we now have two Lyapunov equations with right-hand sides of rank at most
$2 k_d$, which we can exploit when forming the Hadamard products
$\widetilde\DY(\nu)$ and $\widetilde\DW(\nu)$ for the Lyapunov solutions
of~\eqref{eq.lyap2bis}--\eqref{eq.lyap1bis},
\begin{align}
        \DY(\nu) &= T\widetilde\DY(\nu)T^H,
        \qquad\quad
        \widetilde\DY(\nu) =
        L(\overline{\L})\circ
        T^{-1}(U_{\nu}^{}U_\nu^T Y_0^{}+Y_0^{} U_\nu^{}U_{\nu}^{T} )T^{-H}
        \label{solDY}\\
        \DW(\nu) & = T^{-H} \widetilde\DW(\nu)T^{-1},
        \quad
        \widetilde\DW(\nu) =
        L(\L)\circ T^H(U_{\nu}^{}U_\nu^T W_0+W_0^{}U_\nu^{}U_{\nu}^{T})T.
        \label{solDW}
\end{align}
Note that in practice, we can expect the total number of dampers $k_d$ to be
such that $k_d\le s\le n$.

To construct $\DY(\nu)$ and $\DW(\nu)$, we need $T$ and its inverse.
The matrix  $\Comph(\nu)$ is a linearization of
\begin{equation}\label{def.Qt}
        \Qt_\nu(\l )= \Phi^T Q_\nu(\l) \Phi =
        \l^2 I_n +
        \l (\Gamma + R\S_\nu R^T)
        + \Omega^2,
\end{equation}
so both $\Comph(\nu)$ and $\Qt_\nu(\l)$ share the same eigenvalues $\L$.
Also, if $V_R$ is a matrix of right eigenvectors of $\Qt_\nu(\l)$ then
$(\L,V_R)$ is a standard pair for $\Qt_\nu(\l)$ \cite{glr09}, that is,
$\twobyone{V_R}{V_R\L}$ is nonsingular and
$$
        V_R\L^2+ (\Gamma+R\S_\nu^{}R^T)V_R\L +\Omega^2 V_R=0.
$$
Hence,
\begin{equation}\label{eq.T}
     T = \begin{bmatrix}\Omega V_R\\V_R\L\end{bmatrix}
\end{equation}
is nonsingular and it is easy to check that $\Comph(\nu)T = T \L$, i.e.,
$T$ diagonalizes $\Comph(\nu)$.

Since $\Qt_\nu(\l)$ is symmetric, $\overline{V_R}$
is a left eigenvector matrix for $\Qt_\nu(\l)$.
Then letting
$T_{\text{inv}} =
  \begin{bmatrix}
  \big(\L  V_R^T+V_R^T(\Gamma+R\S_\nu^{}R^T)\big)\Omega^{-1}
  &V_R^T
  \end{bmatrix}
$,
we find that $T_{\text{inv}} \Comph(\nu) = \Lambda T_{\text{inv}}$.
Let $\Comph(\nu)$ have $p$ distinct eigenvalues so that
$\L=\diag(\l_1 I_{\ell_1},\ldots,\l_r I_{\ell_r})$.
Then as long as $T_{\text{inv}}$ is of full rank,
$S T^{-1} =  T_{\text{inv}}$ for some nonsingular block diagonal matrix $S$.
Partitioning $V_R$ accordingly to $\L$, i.e.,
$V_R = \begin{bmatrix}V_{R,1} & \ldots& V_{R,r}\end{bmatrix}$,
we have that the $j$th  $\ell_j\times \ell_j$ diagonal block of
$S$ is given by
\begin{equation}\label{def.S}
        S_j = V_{R,j}^T (2\l_j I+\Gamma + R\S_\nu^{}R^T) V_{R,j}.
\end{equation}
Hence
\begin{equation}\label{eq.Tinv}
   T^{-1}=
   S^{-1}\begin{bmatrix}
   \big(\L  V_R^T+V_R^T(\Gamma+R\S_\nu^{}R^T)\big)\Omega^{-1}
    &
   V_R^T
  \end{bmatrix}.
\end{equation}

\subsection{Computation of the Objective Function, its Gradient, and Hessian}
\label{sec.compf}

For a given $\nu$, we compute the eigenvalues $\L$ and the right eigenvector
matrix $V_R$ of $\Qt_\nu(\l)$  in $O(k_dn^2)$ operations using a modification of
Taslman's algorithm~\cite{tasl15}  that includes both internal and external
damping (see \cite{liti26} for details).

\subsubsection{Objective Function}\label{sec.f}
Using \eqref{eq.traceWY}, the decomposition $W(\nu)=W_0+\DW(\nu)$, and
\eqref{solDW}, we have that
$$
f(\nu) = \trace\big(Y(\nu)\big)
=
\frac{1}{2s}\trace\big(G^TW(\nu)G\big)
=
\underbrace{\frac{1}{2s}\trace\big(G^TW_0 G\big)}_{f_0}+
\frac{1}{2s} \trace\big(G^H T^{-H} \widetilde\DW(\nu)T^{-1} G\big).
$$
Now $G^TW_0 G=\Psi_1(1\colon s,1\colon s)+\Psi_3(1\colon s,1\colon s)$ is independent of $\nu$ and
its trace only needs to be computed once.
Compute the $2n\times k_d$ matrices
\begin{equation}\label{def1}
E = V_R^T R,
\quad
\Upsilon_{R}=\begin{bmatrix}\Upsilon_2\\ \Upsilon_3\end{bmatrix}R,
\quad
\Psi_{R}=\begin{bmatrix}\Psi_2\\ \Psi_3\end{bmatrix}R,
\end{equation}
and use $E$ to form the block diagonal matrix $S$ in \eqref{def.S}
and $T^{-1}$ in \eqref{eq.Tinv}.
On using the structure of $U_\nu$ in \eqref{def.Unu},
$T$ in \eqref{eq.T}, and $W_0$ in \eqref{eq.Y0W0}, we compute
$\widetilde\DW(\nu)$ as
\begin{equation}\label{eq-Tu}
\widetilde\DW(\nu) = L(\L)\circ(B+B^{H}),
\quad
B= (\overline{\L E\Sigma_\nu})T_\Psi^{H},
\quad
T_\Psi = T^{H}\Psi_R,
\end{equation}
where the scaling $\L E\Sigma_\nu$ has already been computed for $T^{-1}$.
Let $\mathrm{ind} = [(1:s), (n+1: n+s)]$.
Then
$$
f(\nu) =f_0+ \frac{1}{2}
e^T\Big(
\overline{T^{-1}(:,ind)}
\circ \big(\widetilde\DW(\nu)T^{-1}(:,ind)\big)
\Big)e,
$$
has an overall cost of $O(k_d n^2+sn^2)$, where the
$k_d n^2$ is due to the computation of $\L$ and $V_R$ and the products
of $2n\times n$ and $n\times k_d$ matrices, and the
$sn^2$ is due to the matrix multiplication
$\widetilde\DW(\nu)T^{-1}(:,ind)$.

\subsubsection{Gradient}\label{sec.gradient}
Using the decomposition \eqref{def.W0}, we rewrite $e_i^T\nabla f(\nu)$
in~\eqref{eq.nablaf_i} as
$$
e_i^T\nabla f(\nu)
=
-2\trace\big(
U_i^T ( Y_0W_0+T \widetilde\DY(\nu) \widetilde \DW(\nu) T^{-1} +
 Y_0T^{-H} \widetilde\DW(\nu) T^{-1}+
 T \widetilde\DY(\nu) T^{H}W_0) U_i\big).
$$
We compute $\widetilde\DY(\nu)$ as
\begin{equation}\label{eq-Tp}
\widetilde\DY(\nu) = L(\overline\L)\circ(C+C^H),
\quad
C= (S^{-1}E\Sigma_\nu)T_{\Upsilon}^H, 
\quad
T_{\Upsilon} = {T^{-1}}\Upsilon_R,
\end{equation}
where the scaling $S^{-1}E\Sigma_\nu$ is available from the computation
of $T^{-1}$.
Compute the $2n\times k_d$ matrices
$$
E_Y = \widetilde{\DY}(\nu)\L E,
\quad
E_W = \widetilde{\DW}(\nu) S^{-1}E,
$$
and partition their columns as well as those of   $T_{\Psi}$, $T_{\Upsilon}$,  $\Upsilon_R$ and $\Psi_R$ in \eqref{eq-Tu}, \eqref{eq-Tp} and
\eqref{def1}  into $k$ $2n\times r_i$ blocks
$ E_{Yi},  E_{Wi}, T_{\Psi i}, T_{\Upsilon i}, {\Upsilon}_{Ri}, {\Psi}_{Ri}$ , $i=1\colon k$.
Then
$$
e_i^T\nabla f(\nu)= -2e^T(\Upsilon_{Ri}\circ \Psi_{Ri})e
-2e^T(
E_{Yi}\circ E_{Wi}+T_{\Upsilon i}\circ E_{Wi}
+E_{Yi}\circ T_{\Psi i})e.
$$
The first term is independent of $\nu$ and only needs to be computed once.
The computation of $\nabla f(\nu)$ is dominated by matrix products between $2n\times 2n$ and
$2n\times k_d$ matrices so the overall cost is $O(k_d n^2)$.

\subsubsection{Hessian}
The output $\nu_*$ of Algorithm~\ref{alg.BBRMA} or Algorithm~\ref{alg.SPG}
is a strict local minimum if the submatrix $\nabla^2 f(\nu_*)_{I_*I_*}$ with $I_*=\{i\ |(\nu_{*})_{i}>0\}$ of the $k\times k$
Hessian in \eqref{eq.nabla2f} is positive definite,
which can be verified with a Cholesky factorization assuming we can construct $\nabla^2 f(\nu_*)_{I_*I_*}$ efficiently.
Now the $k(k+1)/2$ entries of the upper part of $\nabla^2 f(\nu)$ can also be
computed in $O(k_dn^2)$ operations. We do not provide all the details, just the
main approach.

Let $Z_j$ be the solution of the Lyapunov equation
\be\label{eq-lya-Z}
\Comph(\nu_*)Z_j+Z_j \Comph(\nu_*)^T = U_jU_j^TY_*+ Y_*U_jU_j^T.
\ee
Then $\vec(Z_j) =
\mP(\nu_*)^{-1}\mathcal M^T e_j$ so that
$
e_i^T\mathcal{N} \mP(\nu_*)^{-1}\mathcal M^T e_j
=
\vec(W_*^{}U_i^{}U_i^T+U_i^{}U_i^TW_*^{})^T\vec(Z_j)
=2\trace\big(U_i^TZ_j W_*U_i\big)
$.
Hence
\begin{align*}
\nabla^2 f(\nu_*)_{ij} &=
e_i^T(\mathcal{N} \mP(\nu)^{-1}\mathcal M^T +
\mathcal{M} \mP(\nu)^{-T}\mathcal N^T) e_j\\
&=
2\trace(U_i^TZ_j W_*U_i+U_j^TZ_i W_*U_j)\\
& = 2\trace(U_i^TZ_j W_0U_i+U_i^TZ_j \Delta W_*U_i+U_j^TZ_i W_0U_j+U_j^TZ_i \Delta W_*U_j)\\
&= 2 e^T\big(
\widetilde{Z_j} T^HU_i \circ T^HW_0U_i +\widetilde{Z_j} T^{H}U_i \circ \widetilde{\Delta W_*}T^{-1}U_i+ \widetilde Z_iT^H U_j \circ T^HW_0 U_j+ \widetilde Z_iT^H U_j \circ\widetilde{\Delta W_*}T^{-1} U_j
\big)e\\
& = 2 e^T\big(
\widetilde{Z_j} \L E_i \circ (T_{\Psi i} + E_{Wi})\big) e + 2e^T\big(\widetilde{Z_i} \L E_j \circ (T_{\Psi j}+  E_{W j})
\big)e.\\
\end{align*}
The last equation is due to  the notation and matrices defined in Sections~\ref{sec.f}--\ref{sec.gradient}
and computed at the last iteration of either
Algorithm~\ref{alg.BBRMA} or Algorithm~\ref{alg.SPG},
and
$$
\widetilde{Z}_j = L(\bar{\L})\circ (H+H^H),
\quad
H = S^{-1}E_i(T_{\Upsilon i} +E_{Yi})^H.
$$
%
%

\section{Numerical Experiments}\label{sec.numexp}

%
%
%
%
%
%

\begin{table}
\caption{List of problems and their parameter values for
\texttt{damp1}, \texttt{damp2}, and \texttt{beam} described at the start of
section~\ref{sec.numexp}.}\label{tab.pbs}
\begin{tabular}{@{}lrcl@{}}
\toprule
Problem's name & $n$ &$\kappa$ & dampers' positions\\
\midrule
\texttt{damp1-a} & 4 & 5 & $\ell_1=2$\\
\texttt{damp1-b} & 20&25 & $\ell_1=2$\\
\texttt{damp1-c} & 20&25 & $(\ell_1,\ell_2) = (2,19)$\\
\texttt{damp2-a} &801&-- & $(\ell_1,\ell_2,\ell_3) = (50,550,220)$\\
\texttt{damp2-b}&1601&-- & $(\ell_1,\ell_2,\ell_3) = (50,950,120)$\\
\texttt{damp2-c}&2001&-- & $(\ell_1,\ell_2,\ell_3) = (850,1950,120)$\\
\texttt{beam-a}&200&-- & $(\ell_1,\ell_2,\ell_3) = (50,100,50)$\\
\texttt{beam-b}&1000&-- & $(\ell_1,\ell_2,\ell_3,\ell_4,\ell_5) = (150,300,500,700,850)$\\
\botrule
\end{tabular}
\end{table}

We consider three types of problems, all with $\Dint$ positive definite so as to ensure they are all stable.

\begin{enumerate}
\item[(i)] \texttt{damp1} corresponds to a one-row mass-spring-damper system with
$n$ masses $m_i$, $i\in[n]$,
$n+1$ springs, all with stiffness $\kappa$,
and $k$ dampers with damping coefficients $\nu_i$, $i\in[k]$ and
damper $i$ positioned on the
$\ell_i$th mass as shown, for example, in
\cite[Fig.~1]{btt13}.
The corresponding matrices are given by
\begin{align}
M &= \diag(m_1,m_2,\ldots,m_n), \nonumber\\
K & = \kappa\tridiag(-1,2,1), \nonumber\\
\Dext(\nu)&=\sum_{i=1}^k\nu_i e_{\ell_i}^{}e_{\ell_i}^T,
\label{eq.Dnu}
\end{align}
and $D_{int}$ as in \eqref{eq-cint1}.
For our numerical experiments, we fix $\alpha$ to $0.01$ and $m_i=i$, $i\in[n]$
but vary the value of $\kappa$, the number $k$ of dampers and their locations,
i.e., the $\ell_i$--see Table~\ref{tab.pbs}.
\item[(ii)] \texttt{damp2} corresponds to a two-row mass-spring-damper system
with $n$ masses $m_i$, $i\in[n]$,
$n+2$ springs with three different stiffness
$\kappa_i$, $i\in[3]$, and $k$ dampers with damping coefficients $\nu_i$, $i\in [k]$
as shown, for example, in \cite[Fig.~1]{jst22}.
The corresponding matrices are given by
\begin{align*}
M &= \diag(m_1,m_2,\ldots,m_n), \quad n=2t+1,\\
K &= \begin{bmatrix}
      \kappa_1\widetilde K&&-\kappa_1 e_t\\
     &\kappa_2\widetilde K&-\kappa_2 e_t\\
    -\kappa_1 e_t^T&-\kappa_2e_t^T &  \kappa_1+\kappa_2+\kappa_3
     \end{bmatrix},
     \quad
     \widetilde{K} = \tridiag(-1,2,-1)\in\R^{t\times t}, \\
\Dint(\nu)&=\Dext(\nu)
           =\nu_1e_{\ell_1}^{}e_{\ell_1}^T
           +\nu_2(e_{\ell_2}^{}-e_{\ell_3+ t}^{})(e_{\ell_2}^{}-e_{\ell_3+t}^{})^T
           +\nu_3 e_{\ell_3}^{}e_{\ell_3}^T,
\end{align*}
and $D_{int}$ as in \eqref{eq-cint1}.
For our numerical experiments we choose the same parameters as in
\cite{jst22}, i.e.,
$\alpha=0.01$ for $\Dint$ in \eqref{eq-cint1}, $\kappa_1 = 100$, $\kappa_2 = 150$, $\kappa_3 = 200$, and
$$
m_i = \left\{
\begin{array}{ll}
0-4 i, & i = 1,\dots, t/2,\\
3i-800, & i = t/2+1, \dots, t,\\
500+i, & i = t +1,\dots, 2t,\\
1800, & i = 2t+1.
\end{array} \right.
$$
We use different values for $n$ and location $\ell_j$ of the three dampers--see Table~\ref{tab.pbs}.

\item[(iii)] \texttt{beam} corresponds to a
 damped slender beam simply supported at both ends
as described in \cite{hmtg08}.
The $\nbyn$ matrices $M$ and $K$ are generated with the
MATLAB toolbox NLEVP~\cite{bhms13} using
\texttt{coeffs = nlevp('damped\_beam',n)},
$\texttt{M = coeffs\{3\}}, \texttt{K = coeffs\{1\}}$.
The internal damping $\Dint$ is as in~\eqref{eq-cint1}
and the external damping $\Dext(\nu)$ is as in~\eqref{eq.Dnu}.
For our numerical experiments, we fix $\alpha$ to $0.2$
and vary the number $k$ of dampers and their locations defined through the
$\ell_i$--see Table~\ref{tab.pbs}.
\end{enumerate}

\subsection{BBRMA versus SPG}

\begin{table}[t]
\caption{A comparison between BBRMA and SPG for the problems in
Table~\ref{tab.pbs}, where $n$ is the size of the problem,
$k$ is the number of dampers, $s$ is the number of eigenmodes to dampen
(the smallest ones),
\#iter is the total number of iterations,
the initial vector $\nu^{(0)}$ is a multiple of the vector of all ones $e$,
{\rm \#ls} is the number of times line search used by SPG,
and {\rm \#eig} is the total number of eigendecompositions.
}\label{table.BBRMAversusSPG}
\begin{tabular}{@{}lrcrclcccc@{}}
\toprule
Problem & $n$ & $k$ & $s$ &$\nu^{(0)}$ & method
&\#iter& \#ls &\#eig & $f(\nu_*)$ \\
\midrule
\texttt{damp1-a} &4&1&4&$e$ & BBRMA & 1000&& 1001&4.3e0\\
          & & & &$e$&SPG   &11& 2& 14 &3.6e0 \\
\midrule
\texttt{damp1-b} &20&1&20&$e$&BBRMA &29&&30 &2.1e1 \\
          &  & &  &$e$&SPG   &10&1&12 &2.1e1 \\
\midrule
\texttt{damp1-c}&20&2&20&$10e$&BBRMA&29&&30 &1.0e1 \\
          &  & &  &$10e$&SPG  & 29&0&30& 1.0e1\\
\midrule
\texttt{damp1-c}&20&2&20&$e$&BBRMA&--&  &-- &-- \\
         &  & &  &$e$&SPG  & 170&43&259& 1.0e1\\
\midrule
\texttt{damp2-a} &801&3&27&$100e$&BBRMA& 24 &&25&1.1e3 \\
       &   & &  &$100e$&SPG  &24&0&25 & 1.1e3\\
\midrule
\texttt{damp2-b}&1601&3&27&$100e$&BBRMA& 28&&29 &3.5e3\\
       &    & &  &$100e$&SPG  & 28&0&29&3.5e3 \\
\midrule
\texttt{damp2-c} &2001&3&20&$100e$&BBRMA&20&&21&3.8e3 \\
       &    & &  &$100e$&SPG  &20&0&21 & 3.8e3\\
\midrule
\texttt{beam-a} &200&3&40&$e$&BBRMA&21&&22& 1e-3\\
       &    & &  &$e$&SPG  &21&0&22 & 1e-3\\
\midrule
\texttt{beam-b} &1000&5&80&$e$&BBRMA&33&&34 & 4e-4 \\
       &    & &  &$e$&SPG  &33&0&34 & 4e-4  \\

\botrule
\end{tabular}
%

\end{table}

For both algorithms, the iteration is terminated when either
the maximum number of iterations $\mathrm{iter}_{\max}$ is reached or
when both criteria
\begin{equation}\label{def.criteria}
\normt{h\big(\nu^{(j)}\big)}<\tolr,
\qquad
\normt{\nu^{(j)}-\nu^{(j-1)}}\le \toln \normt{\nu^{(j-1)}}
\end{equation}
are satisfied.
The first criterion corresponds to line~2 in Algorithms~\ref{alg.BBRMA}
and~\ref{alg.SPG}, while the second tests convergence to an accumulation point.
In all numerical experiments we set
$\mathrm{iter}_{\max}=1000$, $\tolr=10^{-8}$, and $\toln=10^{-5}$.

Table~\ref{table.BBRMAversusSPG} reports the number of iterations required by
each method, together with the number of times nonmonotone line search
invoked by SPG.
Although the eigenvalue decomposition of $\Comph(\nu)$ (computed via that of
$\Qt_\nu(\l)$ in \eqref{def.Qt}) has relatively low complexity,
it involves few BLAS~3 operations and therefore constitute the most
time-consuming component in the evaluation of
$f(\nu)$ and $\nabla f(\nu)$.
For this reason, we use the total number of eigenvalue decompositions as a proxy
for execution time and report it in Table~\ref{table.BBRMAversusSPG}
for both algorithms.
The corresponding optimal damping coefficients are listed in
Table~\ref{table.SPGversusFODA}.

When no line search is required, the total number of eigenvalue decompositions
equals the number of iterations plus one for both BBRMA and SPG, the additional
decomposition being needed to compute $h(\nu^{(0)})$.
This situation occurs for the \texttt{damp2-x} and \texttt{beam-x} test
problems, for which BBRMA and SPG exhibit similar behaviour.

In contrast, for \texttt{damp1-x} problems, SPG frequently activates the
nonmonotone line search. Although this increases the cost per iteration, it can
significantly reduce the total number of iterations and consequently the
overall number of eigenvalue decompositions, as observed for \texttt{damp1-a}
and \texttt{damp1-b}. The choice of the initial damping vector
$\nu^{(0)}$ also has a notable impact on convergence. For instance, in
\texttt{damp1-c}, BBRMA fails to converge when
$\nu^{(0)}=e$
but converges when
$\nu^{(0)}=10e$.
In contrast, SPG converges for both initializations: when
$\nu^{(0)}=e$,
the line search is invoked 43 times and convergence is achieved with a larger
number of iterations than for $\nu^{(0)}=10e$.
Since each line search may require multiple eigenvalue decompositions, this
explains why, for \texttt{damp1-c} with
$\nu^{(0)}=e$,
the total number of eigenvalue decompositions exceeds the sum of the number of
iterations and line-search invocations.

\subsection{SPG versus FODA}
\begin{table}[t]
\caption{A comparison between FODA and SPG for some of the problems in
Table~\ref{tab.pbs}, where the triple $(n,k,s)$ provides
the size of the system, the number of dampers,
and the number of eigenmodes to dampen.
The initial vector $\nu^{(0)}$ is a multiple of the vector of all ones $e$,
\#iter is the total number of iterations,
$\nu_*$ is the vector of optimal damping coefficients,
$\mathrm{res}_{\nabla f}= \norm{h(\nu_*)}$ for SPG and
$\mathrm{res}_{\nabla f}
=\norm{\widetilde{\nabla f}(\nu_*)}$ for FODA,
$res_f=|f(\nu^{(j)})-f(\nu^{(j-1)})|/|f(\nu^{(j-1)})|$
at the last iteration $j$,
and $res_\nu=\|\nu^{(j)}-\nu^{(j-1)}\|$.
}\label{table.SPGversusFODA}
\footnotesize
\begin{tabular}{@{}lcclcccccc@{}}
\toprule
Problem & $(n,k,s)$ &$\nu^{(0)}$ &method &\#eig & $\nu_*$
&$\mathrm{res}_{\nabla f}$ & $res_f$ & $res_\nu$& $f(\nu_*)$\\
\midrule
\texttt{damp1-a} &(4,1,4)& $e$ &FODA&23&[4.4]&6e-6&{\bf 3e-7}&4e-2&3.6e0\\
    && &SPG&12&[4.4]&5e-5&2e-6&{\bf 1e-2}&3.6e0\\
\midrule
\texttt{damp1-b}&(20,1,20)&$e$&FODA&159&[-6.1]&9e1&5e-1&{\bf 6e-3}&-3.8e6\\
&&&SPG&11&[18.9]&2e-7 &{\bf 2e-9}&{\bf 2e-3}&2.1e1\\
\midrule
\texttt{damp1-c} &(20,2,20)&$10e$&FODA&59&[9.6,39.3]&6e-4&{\bf 2e-7}&4e-2&1.0e1\\
&&&SPG&24&[9.6,39.3]&4e-5 &{\bf 6e-8}&{\bf 6e-3}&1.0e1\\
\midrule
\texttt{damp1-c}  &(20,2,20)&$e$&FODA&-&[1,1]&-&-&-&-\\
&&&SPG&254&[9.6,39.3]&4e-6 &{\bf 9e-9}&{\bf 2e-3} & 1.0e1\\
\midrule
\texttt{damp2-a}&(801,3,27)&$100e$&FODA&87&[568,385,284]&3e-8&{\bf 2e-8}&2e-1&1.1e3\\
&&&SPG&14&[565,385,284]&6e-4&{\bf 3e-7}&6e-1&1.1e3\\
\midrule
\texttt{damp2-b}&(1601,3,27)&$100e$&FODA&109&[807,1696,422]&1e-2&{\bf 3e-7}&3e-1&3.5e3\\
&&&SPG&21&[807,1694,422]&4e-5&{\bf 2e-9}&1e-1&3.5e3\\
\midrule
\texttt{damp2-c} &(2001,3,20)&$100e$&FODA&93&[637,704,634]&1e-2&{\bf 3e-9}&2e-1&3.8e3\\
&&&SPG&14&[637,704,663]&5e-4 &{\bf 9e-8}&5e-1 & 3.8e3\\
\botrule
\end{tabular}
\end{table}

We next compare the proposed SPG algorithm with the fast optimal damping
algorithm (FODA) of Jakovčević Stor et al. ~\cite{jst22}.
For this comparison, we use their publicly available Julia
implementation\footnote{
\texttt{https://github.com/ivanslapnicar/FastOptimalDamping.jl.}},
which employs a conjugate gradient method to solve the minimization
problem~\eqref{def.critbb_dual} without imposing a lower bound on the damping
vector $\nu$.
In FODA, the gradient $\nabla f\big(\nu^{(j)}\big)$ is approximated by a finite
difference estimate $\widetilde{\nabla f}\big(\nu^{(j)}\big)$,
and the iterations are terminated when at least one of the following stopping
criteria is satisfied:
\begin{equation}\label{def.criteriaFODA}
\normt{\widetilde{\nabla f}\big(\nu^{(j)}\big)}<\tolr,
\qquad
\frac{|f(\nu^{(j)})-f(\nu^{(j-1)})|}{|f(\nu^{(j-1)})|}\le \tolt,
\qquad
\|\nu^{(j+1)}-\nu^{(j)}\|\le \toln.
\end{equation}
To ensure a comparison that is as fair as possible, and for this experiment
only, we adopt the same stopping criteria for SPG, except for the first
condition in~\eqref{def.criteriaFODA}, which we replace by
$\norm{h\big(\nu^{(j)}\big)}\le \tolr$.
Both criteria are similar in nature but the latter is essential for SPG.
We choose $\tolr = 10^{-8}$,
$\tolt=10^{-6}$, and $\toln=10^{-2}$, and report the results in
Table~\ref{table.SPGversusFODA}.

The numerical results show that SPG consistently requires fewer eigenvalue
decompositions than FODA and is therefore computationally faster. Moreover, FODA
may return negative damping coefficients, as observed for the \texttt{damp1-b}
test problem, and does not always converge, as illustrated by the
\texttt{damp1-c} case with $\nu^{(0)}=e$.
In contrast, for all test problems reported in Table~\ref{table.SPGversusFODA},
we verified that the optimal vector $\nu_*$ computed by SPG corresponds to a
strict local minimizer.


\section{Concluding Remarks}

We have addressed the problem of computing optimal damping coefficients for a
damped vibrational system. In Section~\ref{sec.stable}, we identified conditions
under which the system is either never stable or may lose stability for certain
choices of the damping coefficients $\nu$.
In the latter case, we proposed replacing the constraint $\nu\ge 0$ with $\nu\ge
d$, where $d$ is a nonzero nonnegative vector chosen to ensure stability, and
solving
$$
\min_{\nu\in\R^k\atop \nu\ge d}
        \big\{\trace\big(Y(\nu) \big):
        \ \mbox{$Y(\nu)$ solves~\eqref{eq.lyap2}}
        \big\}
$$
in place of~\eqref{def.critbb_dual}.
This modification does not affect the expressions for the gradient and Hessian
of the objective function
$f(\nu)=\trace\big(Y(\nu)\big)$ derived in Section~\ref{sec-property}.
However, it alters the KKT conditions~\eqref{def.kkt_cond},
leading to a modified residual function
$$
        h(\nu)=(\nu-d)-\hbox{max}\left(\nu-d-\nabla f(\nu),0\right)
$$
in BBRMA and SPG, and a modified search direction in SPG,
$$
        d^{(j)}=\hbox{max}\left(\nu^{(j)}-d-\eta^{(j)}\nabla
        f(\nu^{(j)}),0\right)-\nu^{(j)}.
$$
To minimize this residual function, we proposed two algorithms, BBRMA and SPG.
BBRMA offers a favorable balance between simplicity and computational
efficiency, but it is not guaranteed to converge in all cases. In contrast, SPG
employs a nonmonotone line search and is globally convergent. By exploiting
matrix multiplications and an eigenvalue decomposition of a highly structured
matrix, we showed how to compute $f(\nu)$ in $O((s+k_d) n^2)$ operations
and its gradient in $O(k_d n^2)$ operations, where $k_d\ll n$ is the total number of dampers
in the structure, $s\le n$ is the number of modes to be damped, and $n$ the size
of the problem.
This per-iteration cost is low provided that $n$ is not too large, so that an
initial diagonalization of the mass and stiffness matrices is feasible.

Although SPG was expected to be more computationally expensive than BBRMA  due to the
additional function evaluations required by the nonmonotone line search, our
numerical experiments show that the line search can substantially accelerate
convergence, resulting in a lower overall number of function evaluations.


Finally, we have not investigating
whether or not (i) the objective function $f(\nu)$ is convex,
(ii) its gradient $\nabla f(\nu)$ is Lipschitz continuous, or (iii) the Hessian
$\nabla^2 f(\nu)$ is globally Lipschitz continuous.
Positive results in this direction would provide access to a broader class of
optimization techniques, and we leave these questions for future work.

\section*{Funding}
The first author's research was supported by NSFC No. 12071032 and 12271526
and the China Scholarship Council (CSC).
The second author's research was funded by an EPSRC IS-8 Impact Acceleration
Account (IAA) Funding Stream IAA-IS-8 1019.

\bibliographystyle{plain}
\def\noopsort#1{}\def\l{\char32l}\def\v#1{{\accent20 #1}}
  \let\^^_=\v\def\hbk{hardback}\def\pbk{paperback}

\end{document}